\documentstyle[11pt,leqno]{article}
\newtheorem{theorem}{{\sc Theorem}}
\newcommand{\bt}{\begin{theorem}}
\newcommand{\et}{\end{theorem}}
\newcommand{\newsection}[1]{\setcounter{equation}{0} \setcounter{theorem}{0}
\section{#1}}

\newcommand{\NI}{\noindent}
\newcommand{\bea}{\begin{eqnarray}}
\newcommand{\eea}{\end{eqnarray}}

\def \spec#1 {\mathop{#1}}

\def \b #1 {\bf #1}

\newcommand {\CC}{\centerline}

\newcommand{\ity}{\infty}
\newcommand{\raro}{\rightarrow}

\newcommand{\vsp}{\vskip 1em}

\newcommand{\be}{\begin{equation}}
\newcommand{\ee}{\end{equation}}
\newcommand{\ben}{\begin{eqnarray*}}
\newcommand{\een}{\end{eqnarray*}}

\begin{document}
\CC{\bf{Remarks on Some Conditional Generalized Borel-Cantelli Lemmas}}
\vsp
\CC{\bf{B.L.S. Prakasa Rao}}
\vsp
\CC{\bf{CR Rao Advanced Institute of Mathematics, Statistics and Computer Science}}
\CC{\bf{Hyderabad, India}}
\vsp
\NI{Abstract :} We discuss some conditional generalized Borel-Cantelli lemmas and investigate their quantitative  versions following Arthan and Oliva (arXiv:2012.09942).
\vsp
\NI{Key words :} Borel-Cantelli lemma; Quantitative version.
\vsp
\NI{AMS 2020 Subject Classification :} Primary 60G70.
\vsp
\newsection{Introduction} 

Let $(\Omega, \cal A, P)$ be a probability space and $\{A_n, n\geq 1\}$ be a sequence of events in this probability space. The Borel-Cantelli lemma relates the convergence or divergence of the series $\sum_{n=1}^\ity P(A_n)$ with the probability of the event ``$A_n$ infinitely often" which is the set defined by
$$A_n \; i.o = \cap_{n=1}^\ity\cup_{i=n}^\ity A_i.$$
The event ``$A_n$ infinitely often" can also be represented as the event ``$\limsup A_n$".
\vsp
The classical Borel-Cantelli lemma can be stated in two parts. 
\vsp
\NI{\bf Theorem 1.1} (First Borel-Cantelli Lemma) Let $\{A_n, n\geq 1\}$ be an infinite sequence of events such that $\sum_{n=1}^\ity P(A_n)<\ity.$ Then $P(A_n\; i.o)=0.$
\vsp
\NI{\bf Theorem 1.2} (Second Borel-Cantelli Lemma) Let $\{A_n, n\geq 1\}$ be an infinite sequence of mutually independent events such that $\sum_{n=1}^\ity P(A_n)<\ity.$ Then $P(A_n \;i.o)=1.$
\vsp
Let ${\cal F}$ be a sub-$\sigma$-algebra of $\cal A.$ A sequence of events $\{A_n, n \geq 1\}$ is said to be {\it conditionally independent} given a sub-$\sigma$-algebra ${\cal F}$ if
\be
P(\cap_{i=1}^n A_i|{{\cal F}}) = \prod_{i=1}^n P(A_i|{\cal F}) 
\ee
almost surely for all $n \geq 1.$ 
\vsp
A sequence of random variables $\{X_n, n\geq 1\}$, defined on the probability space $(\Omega, \cal A, P)$, is said to be {\it conditionally independent} given a sub-$\sigma$-algebra  ${\cal F}$ if
\be
P(\cap_{i=1}^n[X_i \leq x_i]|{\cal F}) = \prod_{i=1}^n P(X_i \leq x_i|{\cal F}) 
\ee
almost surely for all $x_i, 1\leq i \leq n, n\geq 1.$
\vsp
If ${\cal F}= \{\phi, \Omega\},$ then the conditional independence reduces to the usual notion of independence of random variables. It is known that independence of a set of events does not imply their conditional independence and conditional independence of a set of events does not imply their independence (cf. Prakasa Rao [1]).
\vsp
Properties of sequences of random variables which are conditionally independent given a sub-$\sigma$-algebra have been studied in Prakasa Rao [1] and Roussas [2]. Yuan and Yang [3], Liu and Prakasa Rao [4], Yuan et al. [5], Liu and Zhang [6], Yuan and Li [7] and Chen and Liu [8]investigated properties of some conditional Borel-Cantelli lemmas and their applications and studied limit theorems for conditionally independent random variables. An interesting example of conditionally independent sequence of random variables is a sequence of exchangeable random variables which become conditionally independent given a suitable sub-$\sigma$-algebra. Prakasa Rao [1] and Roussas [2] give other examples of stochastic models in which conditional independence plays a major role such as non-ergodic models dealing with branching processes (cf. Basawa and Scott [10]). 
\vsp
Proofs of the results presented here are akin to the corresponding results of the Borel-Cantelli lemmas but they are not consequences of those results. Hence they have to be stated separately and proved.
\vsp
Majerek et al. [9] obtained a conditional version of the Borel-Cantelli lemma.
\vsp
\NI{\bf Theorem 1.3} Let $\{A_n, n\geq 1\}$ be a sequence of events and suppose that ${\cal F}$ is a sub-$\sigma$-algebra such that 
$$\sum_{n=1}^\ity P(A_n|{\cal F})<\ity $$
almost surely.Then 
$$P(\limsup A_n|{\cal F}) = 0 $$
almost surely. 
\vsp
The following result is a generalized Borel-Cantelli lemma originally due to Barndorff-Nielsen but with a corrected proof by Balakrishnan and Stepanov [11,12]. Let $A^c$ denote the complement of a set $A.$
\vsp
\NI{\bf Theorem 1.4} Let $\{A_n, n\geq 1\}$ be a sequence of events such that $P(A_n) \raro 0$ as $n \raro \ity.$ If, for some $m\geq 0,$
$$\sum_{n=1}^\ity P(A_n^c A_{n+1}^c \dots A_{n+m-1}^c A_{n+m}) <\ity, $$
then 
$$P(\limsup A_n)=0.$$
\vsp
We now prove a conditional version of this  generalized Borel-Cantelli lemma. 
\vsp
\NI{\bf Theorem 1.5} Let $\{A_n, n\geq 1\}$ be a sequence of events and ${\cal F}$ be a sub-$\sigma$-algebra such that $P(A_n|{{\cal F}}) \raro 0$ almost surely as $n \raro \ity.$ If, for some $m\geq 0,$
$$\sum_{n=1}^\ity P(A_n^c A_{n+1}^c \dots A_{n+m-1}^c A_{n+m}|{\cal F}) <\ity $$
almost surely, then 
$$P(\limsup A_n|{\cal F})=0 $$
almost surely and hence 
$$P(\limsup A_n )=0.$$
\vsp
\NI{\bf Proof:} Note that
\ben
P(\limsup A_n|{\cal F}) &= & P(\cap_{n=1}^\ity \cup_{k=n}^\ity A_k|{\cal F})\\
&=& \lim_{n \raro \ity}P(\cup_{k=n}^\ity A_k|{\cal F})\\
\een
almost surely. However, for any fixed $m>n\geq1,$
\ben
P(\cup_{k=n}^\ity A_k|{\cal F}) &=& P(A_n|{\cal F})+ P(A_n^c A_{n+1}|{\cal F})+ P(A_n^c A_{n+1}^c A_{n+2}|{\cal F})+\dots\\
& \leq & P(A_n|{\cal F})+ P(A_n^c A_{n+1}|{\cal F})+ P(A_n^c A_{n+1}^c A_{n+2}|{\cal F})\\
&&\;\;\;\;+\dots+P(A_n^c\dots A_{n+m-2}^c A_{n+m-1}|{\cal F})\\
&&\;\;\;\;+\sum_{k=n}^\ity P(A_k^c \dots A_{k+m-1}^c A_{k+m}|{\cal F}).\\
\een
The last term given above tends to zero as $n \raro \ity$ since it is the tail sum of the infinite series
$$\sum_{n=1}^\ity P(A_n^c A_{n+1}^c \dots A_{n+m-1}^c A_{n+m}|{\cal F}) <\infty$$
which converges almost surely by hypothesis. 

Furthermore, for a fixed $m \geq 0,$
\bea
\;\;\;\;\\\nonumber
P(A_n|{\cal F})+ P(A_n^c A_{n+1}|{\cal F})+ P(A_n^c A_{n+1}^c A_{n+2}|{\cal F})+\dots + P(A_n^c\dots A_{n+m-2}^c A_{n+m-1}|{\cal F})\\\nonumber
\leq P(A_n|{\cal F})+ P(A_{n+1}|{\cal F})+ P(A_{n+2}|{\cal F})+\dots+P(A_{n+m-1}|{\cal F})
\eea
and the last term on the right side of the inequality tends to zero almost surely by hypothesis. Hence
$$P(\limsup A_n|{\cal F})=0 $$
almost surely which, in turn,  implies that
$$E[P(\limsup A_n|{\cal F})] = P(\limsup A_n)=0.$$
\vsp
\NI{\bf Remarks :} This result is also proved in Lemma 2.4 in Chen and Liu  [8]. 
\vsp
As a consequence of Theorem 1.5, we have the following result.
\vsp
\NI{\bf Theorem 1.6} Let $\{A_n, n\geq 1\}$ be a sequence of events and ${\cal F}$ be a sub-$\sigma$-algebra such that $P(A_n|{\cal F}) \raro 0$ a.s. as $n \raro \ity.$ Further suppose that
$$\sum_{n=1}^\ity P(A_n^c A_{n+1}|{\cal F}) <\ity $$
almost surely. Then 
$$P(\limsup A_n|{\cal F})=0 $$
almost surely and hence 
$$P(\limsup A_n )=0.$$
\vsp
\NI{\bf Proof:} This theorem follows from Theorem 1.5 by choosing $m=1$ in that result.
\vsp
\NI{\bf Theorem 1.7} Let $\{A_n, n\geq 1\}$ be a sequence of events and ${\cal F}$ be a sub-$\sigma$-algebra such that $P(A_n|{\cal F}) \raro 0$ almost  surely as $n \raro \ity.$ Further suppose that
$$\sum_{n=1}^\ity P(A_n A_{n+1}^c|{\cal F}) <\ity $$
almost surely. Then 
$$P(\limsup A_n|{\cal F})=0 $$
almost surely and hence 
$$P(\limsup A_n )=0.$$
\vsp
\NI{\bf  Proof:} This result is a consequence of Theorem 1.6 and the observation that
$$\sum_{n=1}^\ity  P(A_n A_{n+1}^c|{\cal F})= P(A_1|{\cal F})+ \sum_{n=1}^\ity P(A_n^c A_{n+1}|{\cal F}).$$
\vsp
\NI{\bf Remarks:} Result in Theorem 1.7 is  also obtained in Lemma 2.2 in Chen and Liu [8].
\vsp
\newsection{Conditional version of lemma in Balakrishnan and Stepanov [12]}

Balakrishnan and Stepanov [12] generalized the Second result Borel-Cantelli Lemma to some dependent events. Following Balakrishnan and Stepanov [12], given an event $A,$ and a sub-$\sigma$-algebra ${\cal F},$ we say that the quantity $\alpha \geq 0$  is the {\it power-$A$ coefficient} of the conditional independence between the event $A$ and another event $B$ if 
$$P(AB|{\cal F})= [P(A|{\cal F})]^\alpha P(B|{\cal F}) $$
almost surely. It is obvious that if $A$ and $B$ are conditionally independent given ${\cal F},$ then $\alpha=1.$ 
\vsp
Let $\{A_n, n \geq 1\}$ be a sequence of events and let $A^*_n=A_n^c A_{n+1}^c \dots.$ Suppose that the power-$A_n^c$ coefficient of the conditional independence between the events $A_n^c$ and $A^*_{n+1}$ given the sub-$\sigma$-algebra ${\cal F}$ is $\alpha_n.$ Then the following result holds.
\vsp
\NI{\bf Theorem 2.1} Let $A_n$ be a sequence of events as defined above and let ${\cal F}$ be a sub-$\sigma$-algebra . Then 
\be
P(\limsup A_n|{\cal F})=1  
\ee
almost surely if and only if
\be
\sum_{n=1}^\ity \alpha_n P(A_n|{\cal F})= \ity 
\ee
almost surely.
\vsp
\NI{\bf Proof:} Note that
\be
P(A^*_n|{\cal F}) = (P(A_n^c|{\cal F}))^{\alpha_n}P( A^*_{n+1}|{\cal F}) 
\ee
almost surely. Repeating the process, we obtain that
\be
P(A^*_n|{\cal F})= (P(A_n^c|{\cal F}))^{\alpha_n}\dots (P(A_{n+k-1}^c|{\cal F}))^{\alpha_{n+k-1}}P(A^*_{n+k}|{\cal F}) 
\ee
almost surely for all $n \geq 1$ and $k\geq 1.$ Applying the inequality $\log (1-x)\leq -x$ for $0\leq x<1,$ it follows that 

$$ P(A^*_n|{\cal F}) \leq \exp(-\sum_{i=n}^{n+k-1}\alpha_i P(A_i|{\cal F}))P(A^*_{n+k}|{\cal F})  $$
almost surely for all $n \geq 1$ and $k\geq 1.$ Taking limit as $k\raro \ity,$ it follows that
\be
P(A^*_n|{\cal F}) \leq  \exp(-\sum_{i=n}^{\ity}\alpha_i P(A_i|{\cal F})) (1-P(\limsup A_n|{\cal F}))  
\ee
almost surely. Note that (2.2) implies (2.1). Suppose that $P(\limsup A_n|{\cal F}) <1.$ If the series $\sum_{i=n}^{\ity}\alpha_i P(A_i|{\cal F})$ is divergent, then we have a contradiction in (2.5). Hence (2.1) implies (2.2). 
\vsp
Following Balakrishnan and Stepanov [12], let the sequence of events $A_n, n\geq 1$ be called {\it conditionally Markov given a $\sigma$-algebra ${\cal F},$} if the corresponding indicator random variables $I_{A_n}, n \geq 1$ form a Markov chain conditionally given the $\sigma$-algebra ${\cal F}.$ Suppose $\beta_n$ is the power-$A_n^c$ coefficient of conditional dependence of $A_n^c$ and $A_{n+1}^c$ given the $\sigma$-algebra ${\cal F}.$ It can be checked that $\alpha_n$ defined above is equal to $\beta_n$ and in fact 
$$\alpha_n= \frac{\log P(A_n^c A_{n+1}|{\cal F})-\log P(A_{n+1}|{\cal F})}{\log P(A_n^c|{\cal F})}.$$
\vsp
\newsection{Quantitative version of the Conditional Second Borel-Cantelli Lemma}
\vsp
Kochen and Stone [13] presented a result that generalized the Second Borel-Cantelli Lemma  leading to a lower bound on $P(A_n \; i.o)$ when the events $A_n$ are not mutually independent and proved that the Second Borel-Cantelli Lemma holds for pairwise independent. Yan [14] generalized this result leading to the following theorem.
\vsp 
\NI{\bf Theorem 3.1} (Kochen and Stone [13], Yan [14]) Let $\{A_n, n \geq 1\}$ be an infinite sequence of events such that $\sum_{n=1}^\ity P(A_n)=\ity.$ Then
\be
P(A_n \; i.o)\geq \limsup_{n \raro \ity} \frac{[\sum_{k=1}^n P(A_k)]^2}{\sum_{i,k=1}^{n}P(A_iA_k)}.
\ee
\vsp
A related theorem due to Erdos and Renyi  [16] is the following result.
\vsp
\NI{\bf Theorem 3.2} (Erdos and Renyi  [16]) Let $\{A_n, n \geq 1\}$ be an infinite sequence of events such that $\sum_{n=1}^\ity P(A_n)=\ity$ and 
\be
\limsup_{n \raro \ity} \frac{\sum_{i,k=1}^{n}P(A_iA_k)}{[\sum_{k=1}^n P(A_k)]^2}=1.
\ee
Then $P(A_n \; i.o)=1.$
\vsp
We will now obtain the conditional quantitative versions of these results following the ideas and techniques in Arthan and Oliva  [17]. 
\vsp
Given an infinite sequence of events $\{A_n, n \geq 1\}$ and a sub-$\sigma$-algebra ${\cal F},$ the first conditional Borel-Cantelli lemma implies that the probability of the event $\limsup A_n$ conditional on ${\cal F}$ is zero almost surely when $\sum_{n=1}^\ity P(A_n|{\cal F})<\ity$ almost surely. The assumption that $\sum_{n=1}^\ity P(A_n|{\cal F})<\ity$ almost surely is equivalent to the almost sure convergence of the sequence 
$v_k= \sum_{i=1}^k P(A_i|{\cal F}).$ Note that $\{v_k, k\geq 1\}$ is a random sequence. This condition implies the existence of a random function $\eta(.)$ such that
\be
|v_m-v_n|<\frac{1}{2^\ell}, \ell,m,n >\eta(\ell) 
\ee
almost surely since the sequence $\{v_k, k \geq 1\}$ is a Cauchy sequence a.s. In other words, there exists a random function $\phi(\ell)$ such that
\be
\sum_{i=\phi(\ell)}^\ity P(A_i|{\cal F}) \leq \frac{1}{2^\ell} 
\ee
almost surely. We  now state the quantitative version of the conditional First Borel-Cantelli Lemma stated in Theorem 1.3.
\vsp
\NI{\bf Theorem 3.3} Suppose $\{A_n, n \geq 1\} $ is a sequence of events such that the sequence $\sum_{i=1}^m P(A_i|{\cal F})$ converges almost surely with a rate of convergence $\phi(.),$ that is, for all $\ell \geq 0$ and $m >\phi(\ell),$
$$\sum_{i=\phi(\ell)}^m P(A_i|{\cal F})\leq \frac{1}{2^\ell}$$
almost surely. Then the sequence $\{P(\cup_{i=1}^m A_i|{\cal F}), m\geq 1\} $ converges to 0 almost surely with the same rate, that is, for all $\ell \geq 0$ and $m > \phi(\ell),$
$$ P(\cup_{i=\phi(\ell)}^m A_i|{\cal F})\leq \frac{1}{2^\ell} $$
almost surely.
\vsp
\NI{\bf Proof:} This result is an easy consequence of the fact
$$ P(\cup_{i=\phi(\ell)}^m A_i|{\cal F})\leq \sum_{i=\phi(\ell)}^m P(A_i|{\cal F}) a.s.$$
for all $\ell \geq 0$ and $m >\phi(\ell)$ almost surely.
\vsp
The following theorem gives the quantitative version of the second conditional Borel-Cantelli Lemma.
\vsp
\NI{\bf Theorem 3.4} Suppose $\{A_n, n \geq 1\} $ is a sequence of conditionally independent events given a sub-$\sigma$-algebra ${\cal F}.$ Suppose further that the sequence $\{\sum_{i=1}^nP(A_i|{\cal F}), n\geq 1\}$ diverges with rate $\psi(.)$ almost surely, that is, for all $N\geq 1,$
$$\sum_{i=1}^{\psi(N)}P(A_i|{\cal F})\geq N $$ 
almost surely. Then, for all $n \geq 1,$ and $N \geq 1,$
$$P(\cup_{i=1}^{\psi(n+N-1)} A_i|{\cal F})\geq 1-e^{-N}$$
almost surely.
\vsp
\NI{\bf Proof:} We choose a fixed $n\geq 1$ and $N\geq 1.$ Let $A^c$ denote the complement of an event $A.$ By the conditional independence of the events $\{A_n, n \geq 1\}$ given the sub-$\sigma$-algebra ${\cal F},$ it follows that the events $\{ A_n^c, n\geq 1\}$ are also conditionally independent given sub-$\sigma$-algebra ${\cal F},$ and hence
\bea 
P(\cap_{i=n}^{\psi(n+N-1)}A_i^c|{\cal F}) &=& \Pi_{i=n}^{\psi(n+N-1)}P(A_i^c|{\cal F}) \\\nonumber
&=& \Pi_{i=n}^{\psi(n+N-1)}(1-P(A_i|{\cal F})) \\\nonumber
\eea
almost surely. Taking logarithms on both sides of the equation given above, it follows that
\bea
\log (P(\cap_{i=n}^{\psi(n+N-1)}A_i^c|{\cal F})) &= &\log(\Pi_{i=n}^{\psi(n+N-1)}(1-P(A_i|{\cal F})) \\\nonumber
&=& \sum_{i=n}^{\psi(n+N-1)}\log(1-P(A_i|{\cal F}))  \\\nonumber
& \leq & -\sum_{i=n}^{\psi(n+N-1)}P(A_i|{\cal F}) \\\nonumber
& \leq & -N \\\nonumber
\eea
almost surely by the elementary inequality $\log(1+x) \leq x$ for $x \in (-1,\ity).$ Hence
\be
P(\cap_{i=n}^{\psi(n+N-1)}A_i^c|{\cal F})\leq e^{-N} 
\ee
almost surely which implies that
\be
P(\cup_{i=1}^{\psi(n+N-1)} A_i|{\cal F})\geq 1-e^{-N} 
\ee
almost surely.
\vsp
\NI{\bf Remarks} It is easy to check that the Second Conditional Borel-Cantelli lemma is a consequence of Theorem 3.4. This can be seen from the following arguments. Suppose that $\sum_{i=1}^\ity P(A_i|{\cal F})= \ity $ almost surely. Then there exists a function $\psi(.)$ such that, for all,$ N\geq 1,$
$$\sum_{i=1}^{\psi(N)}P(A_i|{\cal F}) \geq N $$
almost surely. Applying Theorem 3.3, it follows that
$$P(\cup_{i=n}^{\psi(n+N-1)} A_i|{\cal F})\geq 1-e^{-N} $$
almost surely which in turn shows that
$$P(\cup_{i=n}^\ity A_i|{\cal F})\geq 1-e^{-N} $$
almost surely. Hence
$$P(\limsup A_n|{\cal F})=1 $$
almost surely. 
\vsp
\newsection{Quantitative version of the Conditional Erdos-Renyi  theorem}

We will  now state and prove a lemma which will be used later.

\NI{\bf Lemma 4.1:} For any sequence of events $\{A_n,n \geq 1\}$ and a sub-$\sigma$-algebra ${\cal F},$ and for all $n \geq 1,$
\be
\frac{\sum_{i,k=1}^n P(A_iA_k|{\cal F})}{(\sum_{k=1}^nP(A_k|{\cal F}))^2}\geq 1 
\ee
almost surely.
\vsp
\NI{\bf Proof:} Let $\alpha_i= P(A_i|{\cal F})$ and $\eta_n=\sum_{i=1}^n\alpha_i.$ It is obvious that  
$$ E(\eta_n^2|{\cal F})\geq (E(\eta_n|{\cal F}))^2 $$
from the elementary property that the conditional variance of any random variable is greater than or equal to zero whenever it exists. Furthermore
$$ E(\eta_n^2|{\cal F})= \sum_{i,k=1}^n P(A_iA_k|{\cal F}) $$
almost surely and 
$$ E(\eta_n|{\cal F})= \sum_{k=1}^nP(A_k|{\cal F}) $$
almost surely. Hence
$$ \sum_{i,k=1}^n P(A_iA_k|{\cal F})\geq (\sum_{k=1}^n P(A_k|{\cal F}))^2 $$
or equivalently
$$ \frac{\sum_{i,k=1}^n P(A_iA_k|{\cal F})}{(\sum_{k=1}^nP(A_k|{\cal F}))^2}\geq 1 $$
almost surely.
\vsp
Following the ideas of Arthan and Oliva [17] again, we can obtain the following quantitative version of the conditional Erdos-Renyi theorem using Lemma 4.1. We omit the details.
\vsp
\NI{\bf Theorem 4.2} Consider a sequence of events $\{A_n,n \geq 1\}$ and a sub-$\sigma$-algebra ${\cal F}.$ Suppose there exists a random function $\psi(n)$ such that, for all $N\geq 1,$
$$ \sum_{i=1}^{\psi(N)} P(A_i|{\cal F})\geq N $$
almost surely and further suppose that there exists a random function $\phi(\ell,n)$ such that for all $\ell, n$ with $\phi(\ell,n)\geq n,$
$$\frac{\sum_{i,k=1}^{\phi(\ell,n)}P(A_iA_k|{\cal F})}{(\sum_{k=1}^{\phi(\ell,n)}(P(A_k|{\cal F}))^2}\leq 1+\frac{1}{2^\ell} $$
almost surely. Define $n_1=\phi(1,1)$ and for $k>1,$ let $n_k= \phi(k, \max(n_{k-1},k)).$ Then, for all $n\geq 1$ and $\ell \geq 1,$
\be
P(\cup_{i=n}^{n_m}A_i|{\cal F})\geq 1-\frac{1}{2^\ell} 
\ee
almost surely where $m=\max(\psi(2n), \ell + 3).$

\newsection{Quantitative version of Conditional Kochen-Stove theorem}

We first state a lemma.

\NI{\bf Lemma 5.1} For any sequence of real numbers $\{a_n, n \geq 1\},$ events $\{A_n, n\geq 1\}$ and a sub-$\sigma$-algebra ${\cal F},$, the following are equivalent.
\be
P(A_n \; i.o|{\cal F})\geq \limsup a_n 
\ee
almost surely and, for every $m\geq 1, \ell \geq 1,$ there exists  $n>m$ such that for every $j>n$
\be
P(\cup_{i=m+1}^n A_i|{\cal F})+\frac{1}{2^\ell}\geq a_j 
\ee
almost surely.
\vsp
Proof of this lemma is along the same lines as the proof of Lemma A.2 in Arthan and Oliva [17].  The following conditional version of Kochen-Stone inequality can be obtained by arguments similar to those in Kochen and Stone [13]. We now state an equivalent quantitative version.
\vsp
\NI{\bf Lemma 5.2} Let $\{A_n, n\geq 1\}$ be an infinite sequence of events and ${\cal F}$ be a sub-$\sigma$-algebra such that $\sum_{i=1}^\ity P (A_i|{\cal F})= \ity$ almost surely. Then, for every $m\geq 1, \ell \geq 1,$ there exists  $n>m$ such that, for every $j>n,$
\be
P(\cup_{i=m+1}^n A_i)+\frac{1}{2^\ell}\geq \frac{(\sum_{k=1}^j P(A_k|{\cal F}))^2}{\sum_{i,k=1}^j P(A_iA_k|{\cal F})} 
\ee 
almost surely.
\vsp 
\NI{\bf Lemma 5.3} (Conditional Chung-Erdos [15]  inequality) Let $\{A_i, 1\leq i \leq n\}$ be an finite sequence of events and $\cal F$ be a sub-$\sigma$-algebra. Then
\be
P(\cup_{k=1}^n A_k|{\cal F}) \geq \frac{(\sum_{k=1}^j P(A_k|{\cal F}))^2}{\sum_{i,k=1}^j P(A_iA_k|{\cal F})} 
\ee
almost surely.
\vsp
Proof of Lemma 5.3 can be given is along the same line as in Chung and Erdos [15] and Yan [14]. 
\vsp
We now state the quantitative version of conditional Kochen-stone theorem. We omit the details. Proof runs along the same lines as the proof of Theorem 4.2 given in Arthan and Oliva [17] in the unconditional case using Lemmas 5.2 and 5.3.
\vsp
\NI{\bf Theorem 5.4} (Quantitative version of Conditional Kochen-Stone theorem) Let $\{A_n, n\geq 1\}$ be an infinite sequence of events and ${\cal F}$ be a sub-$\sigma$-algebra.. Suppose there exists a random function $\phi(n)$ such that, for all $N\geq 1,$
$$\sum_{i=1}^{\phi(N)}P(A_i|{\cal F})\geq N $$
almost surely. Then, for all $m$ and $\ell$ and $g(.)$ such that $g(i) >i$ for all $i$, there exists $n>m$ such that for all $j\in [n,g(n)],$
\be
P(\cup_{k=1}^n A_k|{\cal F}) +\frac{1}{2^\ell} \geq \frac{(\sum_{k=1}^j P(A_k|{\cal F}))^2}{\sum_{i,k=1}^jP(A_iA_k|{\cal F})} 
\ee
almost surely.
\vsp
\NI{\bf Acknowledgment} This work is supported by the scheme ``INSA Senior Scientist" at the CR Rao Advanced Institute of Mathematics, Statistics and Computer Science, Hyderabad, India.
\vsp
\NI{\bf References}
\begin{description} 

\item 1. Prakasa Rao, B.L.S. (2009) Conditional independence, conditional mixing and conditional association, {\it Ann. Inst. Statist. Math.}, {\bf 61}, 441-460.

\item 2. Roussas, G.G. (2008) On conditional independence, mixing and association, {\it Stoch. Anal. Appl.}, {\bf 26}, 1274-1309.

\item 3. Yuan, D. and Yang, Y. (2011) Conditional versions of limit theorems for conditionally associated random variables, {\it J. Math. Anal. Appl.}, {\bf 376}, 282-293.

\item 4. Liu, J. and Prakasa Rao, B.L.S. (2013) On conditional Borel-Cantelli lemmas for sequences of random variables, {\it J. Math. Anal. Appl.}, {\bf 399}, 156-165.

\item 5. Yuan, D. Wei, L., and Lei, L. (2014) Conditional central limit theorems for a sequence of conditionally independent random variables, {\it J. Korean Math. Soc.}, {\bf 51}, 1-15.

\item 6. Liu, J.C. and Zhang, L.D. (2014) Conditional Borel-Cantelli lemmas and conditional strong law of large numbers, {\it Acta Math. 
Appl. Sin.}, {\bf 37}, 537-546. 

\item 7. Yuan, D. and Li, S. (2015) From conditional independence to conditionally negative association: some preliminary results., {\it Comm. Stat.. Theory and Methods}, {\bf 44}, 3942-3966.

\item 8. Chen, Q and Liu, J. (2017) The conditional Borel-Cantelli lemma and applications, {\it J. Korean Math. Soc.}, {\bf 54}, 441-460.

\item 9. Majerak, D., Nowak, W., and Zieba, W. (2005) Conditional strong law of large numbers,  {\it Int. J. Pure Appl. Math.}, {\bf 20}, 143-156.

\item 10. Basawa, I.V.  and Scott, D. (1983) {\it Asymptotic Optimal Inference for Non-ergodic Models}, Lecture Notes Stat. Vol. 17, Springer, New York. 

\item 11. Balakrishnan, N. and Stepanov, A. (2010) Generalization of Borel-Cantelli lemma, {\it The Mathematical Scientist}, {\bf 35}, 61-62.

\item 12. Balakrishnan, N. and Stepanov, A. (2021) A note on the Borel-Cantelli lemma, arXiv: 2112.00741v1 [math.PR] 1 Dec 2021.

\item 13. Kochen, S.B. and Stone, C.J. (1964) A note on the Borel-Cantelli lemma, {\it Illinois J. Math.} {\bf 8}, 248-251.

\item 14. Yan, J. (2006) A simple proof of two generlized Borel-Cantelli lemmas, In {\it Lecture Notes in Math.}, Vol. 1874, pp. 77-79, Springer.

\item 15. Chung, K.L. and Erdos, P. (1952) On the application of the Borel-Cantelli lemma, {\it Trans. Amer. Math. Soc.}, {\bf 72}, 179-186.

\item 16. Erdos. P. and Renyi, A. (1959) On Cantor's series with convergence $\sum 1/q_n.$, {\it Ann. Univ. Sci. Budap. Rolando Eotvos, Sect. Math.}, {\bf 2}, 93-109.

\item 17. Arthan, R. and Oliva, P. (2020) On the Borel-Cantelli lemmas, the Erdos-Renyi theorem and the Kochen-stone theorem, arXiv:201209942v1 [math. PR] 17 Dec 2020.

\end{description}

\end{document}